\documentclass[12pt]{amsart}
\usepackage{a4wide}
\usepackage{amssymb}
\usepackage{amsthm}
\usepackage{amsmath}
\usepackage{amscd} \usepackage{verbatim}
\usepackage[all]{xy}

\textheight8.3in \textwidth6.3in \numberwithin{equation}{section}
\theoremstyle{plain}
\newtheorem{theorem}{Theorem}[section]
\newtheorem{corollary}[theorem]{Corollary}
\newtheorem{lemma}[theorem]{Lemma}
\newtheorem{proposition}[theorem]{Proposition}

\theoremstyle{definition}

\newtheorem*{example}{Example}

\theoremstyle{remark}
\newtheorem{remark}{Remark}

\newcommand{\ad}{a_{\Delta}}
\newcommand{\SL}{\text {\rm SL}}
\newcommand{\R}{\mathbb{R}}
\newcommand{\Q}{\mathbb{Q}}
\newcommand{\Z}{\mathbb{Z}}
\newcommand{\N}{\mathbb{N}}
\newcommand{\C}{\mathbb{C}}
\newcommand{\ord}{\text {\rm ord}}

\renewcommand{\H}{\mathbb{H}}

\newcommand{\zxz}[4]{\begin{pmatrix} #1 & #2 \\ #3 & #4 \end{pmatrix}}

\newcommand{\leg}[2]{\left( \frac{#1}{#2} \right)}
\newcommand{\kzxz}[4]{\left(\begin{smallmatrix} #1 & #2 \\ #3 & #4\end{smallmatrix}\right) }
\newcommand{\kabcd}{\kzxz{a}{b}{c}{d}}

\newcommand{\PP}{\mathbb{P}}

\newcommand{\calF}{\mathcal{F}}

\newcommand{\calM}{\mathcal{M}}

\newcommand{\eps}{\varepsilon}
\newcommand{\bs}{\backslash}

\newcommand{\sgn}{\operatorname{sgn}}

\newcommand{\Sl}{\operatorname{SL}}

\newcommand{\Aut}{\operatorname{Aut}}

\begin{document}

\title[Harmonic weak Maass forms
and the vanishing of Hecke eigenvalues
]{Differential operators for harmonic weak Maass forms
and the vanishing of Hecke eigenvalues}

\author{Jan H. Bruinier, Ken Ono and Robert C. Rhoades}

\address{Fachbereich Mathematik,
Technische Universit\"at Darmstadt, Schlossgartenstrasse 7, D--64289
Darmstadt, Germany} \email{bruinier@mathematik.tu-darmstadt.de}

\address{Department of Mathematics, University of Wisconsin,
Madison, Wisconsin 53706 USA} \email{ono@math.wisc.edu}
\email{rhoades@math.wisc.edu}
\thanks{The second author thanks
the generous support of the National
Science Foundation, and the Manasse family. The third author is grateful
for the support of  a National Physical Sciences Consortium Graduate
Fellowship, and a National Science Foundation Graduate Fellowship.}

\begin{abstract}
For integers $k\geq 2$, we study two differential operators on
harmonic weak Maass forms of weight $2-k$. The operator $\xi_{2-k}$
(resp. $D^{k-1}$) defines a map to the space of weight $k$ cusp
forms (resp. weakly holomorphic modular forms). We leverage these
operators to study coefficients of harmonic weak Maass forms.
Although generic harmonic weak Maass forms are expected to
have transcendental
coefficients, we show that those forms which are ``dual'' under
$\xi_{2-k}$
to newforms with vanishing Hecke eigenvalues (such as CM forms)
have algebraic coefficients.
Using regularized
inner products, we also characterize the image of $D^{k-1}$.
\end{abstract}

\maketitle

\section{Introduction and Statement of Results}

Let $M_k^{!}(\Gamma_0(N),\chi)$ denote the space of integer weight $k$ weakly
holomorphic modular forms on $\Gamma_0(N)$ with Nebentypus
$\chi$. Recall that a {\it weakly
holomorphic modular form} is any meromorphic modular form whose poles
(if any) are supported at cusps.
Weakly holomorphic modular forms naturally sit in spaces
of harmonic weak Maass forms (see Section~\ref{hwmf} for definitions),
more general automorphic forms
which have been a
source of recent interest  due to their connection to Ramanujan's
mock theta functions, Borcherds products, derivatives of
modular $L$-functions, and traces of singular moduli
(see \cite{ BO1, BO2, BOPNAS, BOR, Br, BF, BrO, ZagierBourbaki, Z2}).

In view of these applications, it is natural to investigate
the arithmeticity of the Fourier coefficients of such Maass forms,
and to also
investigate their nontrivial interplay
with holomorphic and weakly holomorphic
modular forms.
In the works above,
one such nontrivial relationship (see Prop. 3.2 of \cite{BF}), involving the
differential operator
\begin{displaymath}
\xi_w:=2i y^w\cdot\overline{\frac{\partial}{\partial \overline{z}}},
\end{displaymath}
plays a central role. It is the fact that
\begin{equation}\label{ximap}
\xi_{2-k}: H_{2-k}(\Gamma_0(N),\chi)\longrightarrow
S_k(\Gamma_0(N), \overline{\chi}).
\end{equation}
Here $H_w(\Gamma_0(N),\chi)$ denotes the space of weight $w$
harmonic weak Maass forms on $\Gamma_0(N)$ with Nebentypus $\chi$,
and $S_w(\Gamma_0(N), \chi)$ denotes the subspace of cusp forms.  It
is not difficult to make this more precise using Fourier expansions.
In particular, every weight $2-k$ harmonic weak Maass form $f(z)$
has a Fourier expansion of the form
\begin{equation}\label{fourier}
f(z)=\sum_{n\gg -\infty} c_f^+(n) q^n + \sum_{n<0} c_f^-(n)
\Gamma(k-1,4\pi |n|y) q^n,
\end{equation}
where $\Gamma(a,x)$ is the incomplete Gamma-function,
$z=x+iy\in \H$, with $x, y\in \R$, and $q:=e^{2\pi i z}$.
A straightforward calculation shows that
$\xi_{2-k}(f)$ has the Fourier expansion
\begin{equation}\label{xiprop}
\xi_{2-k}(f)=-(4\pi)^{k-1}\sum_{n=1}^{\infty}\overline{c_f^{-}(-n)}n^{k-1}q^n.
\end{equation}

As (\ref{fourier}) reveals, $f(z)$
naturally decomposes into two summands
\begin{eqnarray}
\label{fourierhh}
&f^{+}(z):=\sum_{n\gg -\infty} c_f^+(n) q^n,\\
\label{fouriernh}
&f^{-}(z):=\sum_{\substack{n<0}} c_f^-(n)\Gamma(k-1,4\pi|n|y)q^n.
\end{eqnarray}
Therefore,
$\xi_{2-k}(f)$
is given simply in terms of $f^{-}(z)$, the {\it non-holomorphic part}
of $f$.

Here we show that $f^{+}(z)$, the
{\it holomorphic part} of $f$,
is also intimately related to  weakly holomorphic modular forms.
We require the differential operator
\begin{equation}\label{D}
D:=\frac{1}{2\pi i}\cdot \frac{d}{dz}.
\end{equation}

\begin{theorem}\label{thm1}
If $2\leq k\in \Z$ and $f\in H_{2-k}(\Gamma_0(N),\chi)$, then
\begin{displaymath}
D^{k-1}(f)\in M_k^{!}(\Gamma_0(N),\chi).
\end{displaymath}
Moreover, assuming the notation in (\ref{fourier}), we have
\begin{displaymath}
D^{k-1}f=D^{k-1}f^{+}=\sum_{n\gg -\infty} c_f^{+}(n)n^{k-1}q^n.
\end{displaymath}
\end{theorem}

\begin{remark} Theorem~\ref{thm1} is related to classical results
on weakly holomorphic modular forms and Eichler integrals.
Theorem~\ref{thm1} is a generalization of the classical result
on weakly holomorphic modular forms to the context of harmonic
weak Maass forms.
\end{remark}

Theorem~\ref{thm1} implies that
the coefficients $c_f^{+}(n)$, for non-zero $n$,
are obtained by dividing the $n$th coefficient
of some fixed weakly holomorphic modular form by $n^{k-1}$.
Therefore
we are compelled to determine the image of the
map
\begin{displaymath}
D^{k-1}: H_{2-k}(\Gamma_0(N),\chi)\longrightarrow
M_k^{!}(\Gamma_0(N),\chi).
\end{displaymath}
It is not difficult to see that this map is
not generally surjective. Our next result determines the image
of $D^{k-1}$ in terms of ``regularized'' inner products (see Section
\ref{innerprod}).

\begin{theorem}\label{thm2}
If $2\leq k\in \Z$, then the image of the map
\begin{displaymath}
D^{k-1}: H_{2-k}(\Gamma_0(N),\chi) \longrightarrow
M_k^{!}(\Gamma_0(N),\chi)
\end{displaymath}
consists of those
$h\in M_k^{!}(\Gamma_0(N),\chi)$ which are orthogonal to cusp forms
(see Section~\ref{innerprod}) which also have constant term 0 at all cusps of
$\Gamma_0(N)$.
\end{theorem}

Although these results for $D^{k-1}$ suggest that one has complete
information concerning the Fourier coefficients of $f^{+}$, it turns
out that some of the most basic questions remain open. Here we
consider algebraicity. Despite the fact that we have a fairly
complete theory of algebraicity for forms in
$M_k^{!}(\Gamma_0(N),\chi)$, thanks to the $q$-expansion principle,
the theory of Eisenstein series and newforms, this question remains
open for harmonic weak Maass forms. In view of the theory of
newforms, it is natural to restrict our attention to those $f\in
H_{2-k}(\Gamma_0(N),\chi)$ for which $\xi_{2-k}(f)\in
S_k(\Gamma_0(N),\overline{\chi})$ is a Hecke eigenform. For reasons
which will become apparent, we shall concentrate on those forms for
which
\begin{equation}\label{normcondition}
\xi_{2-k}(f)=\frac{g}{\| g \|^2},
\end{equation}
where $g$ is a normalized newform and $\| g\|$ denotes its usual Petersson norm.

To illustrate the nature of this problem, we consider two examples
of Maass-Poincar\'e series which are not weakly holomorphic modular
forms. The Maass-Poincar\'e series (see Section~\ref{Poincare})
$f:=Q(-1,12,1;z)\in H_{-10}(\SL_2(\Z))$ (note. If a Nebentypus
character is not indicated, then it is assumed to be trivial)
satisfies (\ref{normcondition}) for $g=\Delta(z)$, the unique
normalized weight 12 cusp form on the full modular group. The first
few coefficients of its holomorphic part are
\begin{displaymath}
Q^{+}(-1,12,1;z)
\sim q^{-1}-0.04629-1842.89472q-23274.07545q^2-225028.75877q^3-\cdots.
\end{displaymath}
There is little reason to believe  that these coefficients are
rational or algebraic. On the other hand, we shall prove that the
Maass-Poincar\'e series $Q(-1,4,9;z)\in H_{-2}(\Gamma_0(9))$ has the
property that $Q^{+}(-1,4,9;z)$ has rational coefficients. Its first
few terms are
\begin{equation}\label{goodpoincare}
Q^{+}(-1,4,9;z)=q^{-1}-\frac{1}{4}q^2+\frac{49}{125}q^5-\frac{48}{512}q^8-\frac
{771}{1331}q^{11}+\cdots,
\end{equation}
and $f:=Q(-1,4,9;z)$ satisfies (\ref{normcondition}) for
the unique normalized newform in $S_4(\Gamma_0(9))$.

Our next result explains the distinction between these two cases. To
make this precise, let $g\in S_k(\Gamma_0(N),\overline{\chi})$ be a
normalized newform, and let $F_g$ be the number field obtained by
adjoining the coefficients of $g$ to $\Q$. We say that a harmonic
weak Maass form $f\in H_{2-k}(\Gamma_0(N),\chi)$ is {\it good for}
$g$ if it satisfies the following properties (Section~\ref{hwmf} for
definitions):
\begin{enumerate}
\item[(i)]
The principal part of $f$ at the cusp $\infty$ belongs to  $F_g[q^{-1}]$.
\item[(ii)]
The principal parts of $f$ at the other cusps of $\Gamma_0(N)$ are constant.
\item[(iii)]
We have that $\xi_{2-k}(f)=\|g\|^{-2} g$.
\end{enumerate}

\begin{remark} For every such $g$,
Proposition~\ref{goodmaass} will show that there  is
an $f$ which is good for $g$.
Moreover, such an $f$ is unique
up to a weakly holomorphic form
in $M_{2-k}^!(\Gamma_0(N),\chi)$ with coefficients in $F_g$.
Such $f$ can be constructed explicitly using Poincar\'e series (for example,
see Section~\ref{Poincare} for even $k\geq 2$ and trivial Nebentypus).
\end{remark}

\begin{theorem}\label{thm3}
Let $g\in S_k(\Gamma_0(N),\overline{\chi})$ be a normalized newform
with complex multiplication. If $f\in H_{2-k}(\Gamma_0(N),\chi)$ is
good for $g$, then all coefficients of $f^{+}$ are in
$F_g(\zeta_M)$, where $\zeta_M := e^{2\pi i /M}$, and $M = ND$ where
$D$ is the discriminant of the field of complex multiplication.
\end{theorem}

\begin{remark} $ \ \ \ $

\noindent
i) The rationality of $Q^{+}(-1,4,9;z)$ in (\ref{goodpoincare})
is an example of Theorem~\ref{thm3}.
In this case $Q(-1,4,9;z)$ is good for the unique CM newform
in $S_4(\Gamma_0(9))$. We shall discuss this example in detail
in the last section.

\smallskip

\noindent
ii) The field $F_g$ in Theorem~\ref{thm3} is
explicit (see the discussion in Section~\ref{proof1.3}).

\smallskip
\noindent iii) Suppose that $g\in S_k(\Gamma_0(N),\overline{\chi})$
is a normalized newform. If $f\in H_{2-k}(\Gamma_0(N),\chi)$ is good
for $g$, then the proof of Theorem~\ref{thm3} implies that all of
the coefficients of $f^{+}$ belong to $F_g(c^+_f(1))$. It would be
interesting to describe this field in terms of intrinsic invariants
associated to $g$.

\smallskip
\noindent
iv) It is interesting to compare Theorem~\ref{thm3},
which concerns integer weights $2-k$, with the
results in \cite{BrO} which pertain to weight 1/2 harmonic
weak Maass forms.
The first two authors proved that if $g$ is a newform of weight $3/2$
which is orthogonal to all elementary theta series, and if $f$ is
defined analogously as above, then
\begin{displaymath}
\begin{split}
\# \{ n \in \N \
: \ c_f^{+}(n) \ {\text {\rm transcendental}}\}&=+\infty,\\
\# \{ n \in \N \ : \ c_f^{+}(n) \ {\text {\rm algebraic}}\}&=+\infty.
\end{split}
\end{displaymath}
In fact, estimates are obtained for these quantities.
These results are related to the vanishing of derivatives
of quadratic twists of weight 2 modular $L$-functions at $s=1$.

\smallskip
\noindent v) It would be interesting to find an explicit
construction of good harmonic weak Maass forms for CM newforms.
Perhaps there is a construction which is analogous to the case of
the mock theta functions \cite{BO1, BO2, BOR, ZagierBourbaki, Z2}.

\smallskip
\noindent
vi) In the examples we know, it turns out that
the coefficients of $f^{+}$ are actually contained in
$F_g$. It seems natural to ask whether this is true in general.
\end{remark}

The proof of Theorem~\ref{thm3} relies on the fact that some Hecke
eigenvalues of $g$ vanish. A simple generalization of the proof of
Theorem~\ref{thm3} can be used to detect the vanishing of the
Fourier coefficients of a newform.
\begin{theorem}\label{thm:heckeVanishing}
Suppose that $g=\sum_{n=1}^{\infty} c_g(n) q^n \in
S_k(\Gamma_0(N),\overline{\chi})$ is a normalized newform, and
suppose that $f\in H_{2-k}(\Gamma_0(N),\chi)$ is good for $g$. If
$p\nmid N$ is a prime for which $c_g(p)=0$, then
$c_f^{+}(n)$ is algebraic when $\ord_p(n)$ is odd.
\end{theorem}

\begin{remark}
The proof of Theorem~\ref{thm:heckeVanishing} shows that the
coefficients of $f^{+}$ are in an explicit abelian extension of
$F_g$ when $c_g(p)=0$. It seems possible that the coefficients of
$f^{+}$ are always in $F_g$ when there are any vanishing Hecke
eigenvalues. As the next example will show, this is the case when
$N=1$.
\end{remark}

\begin{example}
Here we consider Lehmer's Conjecture on the nonvanishing of
Ramanujan's $\tau$-function, where
\begin{displaymath}
\Delta(z)=\sum_{n=1}^{\infty}\tau(n)q^n.
\end{displaymath}
This example generalizes easily to all level 1 Hecke eigenforms.

Although Theorem~\ref{thm:heckeVanishing} relates Lehmer's
Conjecture to the alleged transcendence of the coefficients, say
$\ad(n)$, of $Q^{+}(-1,12,1;z)$, it turns out that much more is
true. Lehmer's Conjecture is implied by the mere irrationality of
any these coefficients.

We make use of explicit formulas. Using the classical Eisenstein
series $E_4$ and $E_6$ and the classical $j$-function $j(z)$, we
define polynomials $J_m(x)$ by
\begin{equation}
\sum_{m=0}^{\infty}J_m(x)q^m:=
\frac{E_4(z)^2E_6(z)}{\Delta(z)}
\cdot \frac{1}{j(z)-x}=1+(x-744)q+\cdots.
\end{equation}
For each $m$ we then let $j_m(z)=J_m(j(z))$.
If $p$ is prime,
then define the modular functions
\begin{equation}\label{A}
\begin{split}
A_p(z):= \frac{24}{B_{12}}(1+p^{11})
 +j_p(z)-264\sum_{m=1}^{p}\sigma_9(m)j_{p-m}(z),
\end{split}
\end{equation}
\begin{equation}\label{B}
\begin{split}
B_p(z)&:=-\tau(p)\left(
-264+\frac{24}{B_{12}}+j_1(z)\right).
\end{split}
\end{equation}
Here $B_{12}=-691/2730$ is the 12th Bernoulli number, and
$\sigma_9(n):=\sum_{d\mid n}d^9$. Using the principal part of
$Q(-1,12,1;z)$ combined with the fact that $\Delta(z)$ is an
eigenform of the Hecke algebra, one can show (for example, see
\cite{mockdelta}), for primes $p$, that
\begin{displaymath}
\begin{split}
\sum_{n=-p}^{\infty}\left(p^{11}\ad(pn)-\tau(p)\ad(n)+\ad(n/p)\right)q^n
=\frac{A_p(z)+B_p(z)}{E_4(z)E_6(z)}.
\end{split}
\end{displaymath}
These weight $-10$ modular forms have integer coefficients. Now
suppose that $\tau(p)=0$ for a prime $p$. Then $\ad(np)$ is rational
for every $n$ coprime to $p$. Under this assumption, the proof of
Theorem~\ref{thm3} then implies that $\ad(n)$ is rational when
$\ord_p(n)$ is odd.
\end{example}

Due to Theorems~\ref{thm3} and \ref{thm:heckeVanishing}, it is
natural to consider the arithmetic properties of harmonic weak Maass
forms. For brevity, we will be content  with the following result
for certain forms with prime power level and trivial Nebentypus.

\begin{theorem}\label{serrepadic}
Suppose that $p$ is prime, and that $f(z)\in H_{2-k}(\Gamma_0(p^t))$
is good for a newform $g\in S_k(\Gamma_0(p^t))$ with complex
multiplication. If we let
\begin{displaymath}
a:=\min\{d\geq 0\ : \ c_f^{+}(p^dn)=0 \ {\text {\rm for all}}\
n<0 \},
\end{displaymath}
then the following are true:

\smallskip
\noindent
1) The formal $q$-series
\begin{displaymath}
f^*:=\sum_{n=0}^{\infty} c_f^{+}(p^a n)n^{k-1}q^n
\end{displaymath}
is a $p$-adic modular form on $\SL_2(\Z)$ of weight $k$.

\smallskip
\noindent 2) For every positive integer $b$, we have that
\begin{displaymath}
\lim_{X\rightarrow +\infty} \frac{ \# \{ n\leq X \ : \ c_f^{+}(p^a n)n^{k-1}
\not \equiv 0\pmod{p^b}\}}{X}=0.
\end{displaymath}

\smallskip
\noindent 3) If  $p\leq 7$, or  $p\geq 11$ and $k\equiv 4, 6, 8, 10,
14\pmod{p-1}$, then as $p$-adic numbers we have
\begin{displaymath}
\lim_{n\rightarrow +\infty} c_f^{+}(p^{a+n})p^{n(k-1)}=0.
\end{displaymath}
\end{theorem}

\begin{remark}
Theorem~\ref{serrepadic} (2) says that ``almost every'' $c_f^{+}(p^a
n)$ is a multiple of $p^b$. Theorem~\ref{serrepadic} (3) is not a
trivial statement since the coefficients $c_f^{+}(p^{a+n})$ tend to
have unbounded denominators involving increasing powers of $p$.
\end{remark}

In Section~\ref{hwmf} we recall definitions and facts about harmonic
weak Maass forms and their behavior under certain differential operators.
In Section~\ref{thm1proof} we prove Theorem~\ref{thm1}.
In Section~\ref{innerprod} we recall facts about the regularized inner
product, which generalizes Petersson's inner product, and
we prove Theorem~\ref{thm2}. In Section~\ref{proof1.3} we prove
Theorems~\ref{thm3}, \ref{thm:heckeVanishing}
and \ref{serrepadic}, and in Section~\ref{Poincare}
we illustrate Theorems~\ref{thm1} and
~\ref{thm2} using Poincar\'e series.
In the last section, we examine example (\ref{goodpoincare})
in the context of all of the results
above.

\section*{Acknowledgements}
The authors thank the referee for several helpful suggestions and
corrections.

\section{Harmonic weak Maass forms}
\label{hwmf}
Here we recall definitions and facts about harmonic weak Maass
forms.
Throughout, let $z=x+iy\in \H$, the upper-half of the complex plane,
with $x, y\in \R$.
Also, throughout suppose that
$k\in \N$. We define the weight $k$
hyperbolic Laplacian by
\begin{equation}\label{deflap}
\Delta_k := -y^2\left( \frac{\partial^2}{\partial x^2}+
\frac{\partial^2}{\partial y^2}\right) + iky\left(
\frac{\partial}{\partial x}+i \frac{\partial}{\partial y}\right).
\end{equation}
Suppose that $\chi$ is a Dirichlet character modulo $N$. Then
a {\it harmonic weak Maass form of weight $k$ on
$\Gamma_0(N)$} with Nebentypus $\chi$
is any smooth function on $\H$
satisfying:
\begin{enumerate}
\item[(i)]
$f\left(\frac{az+b}{cz+d}\right) = \chi(d)(cz+d)^kf(z)$ for all $
\left(\begin{matrix}a&b\\c&d\end{matrix}\right)\in \Gamma_0(N)$;
\item[(ii)] $\Delta_k f =0 $;
\item[(iii)]
There is a polynomial $P_f=\sum_{n\leq 0} c_f^+(n)q^n \in \C[q^{-1}]$
such that $f(z)-P_f(z) = O(e^{-\eps y})$ as $y\to\infty$ for
some $\eps>0$. Analogous conditions are required at all
cusps.
\end{enumerate}

The polynomial $P_f\in \C[q^{-1}]$ is called the {\em principal part} of $f$ at the corresponding cusp.
We denote the vector space of these harmonic weak Maass forms by $H_k(\Gamma_0(N),\chi)$

\begin{remark}
Note that our definition slightly differs from the one of \cite{BF},
since we assume that the singularities of $f$ at the cusps are supported on
the holomorphic parts of the corresponding Fourier expansions. This
space is denoted by $H^+_k$ in \cite{BF}.
\end{remark}

Recall the Maass raising and lowering operators
(see \cite{BF, Bump}) $R_k$ and $L_k$ on functions $f:\H\to \C$  which are defined by
\begin{align*}
R_k  &=2i\frac{\partial}{\partial z} + k y^{-1} = i \left(\frac{\partial}{\partial x} - i
\frac{\partial}{\partial y}\right)+ ky^{-1},\\
L_k  &= -2i y^2 \frac{\partial}{\partial\bar{z}} = -i y^2
\left(\frac{\partial}{\partial x} + i \frac{\partial}{\partial
y}\right).
\end{align*}
With respect to the Petersson slash operator (see (\ref{slash})),  these
operators satisfy the intertwining properties
\begin{align*}
R_k (f\mid_k\gamma)  &=  (R_k f)\mid_{k+2} \gamma,\\
L_k (f\mid_k\gamma)  &=  (L_k f)\mid_{k-2} \gamma,
\end{align*}
for any $\gamma\in \Sl_2(\R)$.
The Laplacian $\Delta_k$ can be expressed in terms of $R_k$ and $L_k$ by
\begin{equation}
\label{deltalr}
-\Delta_k = L_{k+2} R_k +k = R_{k-2} L_{k}.
\end{equation}
If $f$ is an eigenfunction of $\Delta_k$ satisfying $\Delta_k f=\lambda f$, then
\begin{align}
\label{rev}
\Delta_{k+2} R_k f  &=  (\lambda+k) R_k f,\\
\label{lev}
\Delta_{k-2} L_k f  &=  (\lambda-k+2) L_k f.
\end{align}
For any positive integer $n$ we put
\[
R_k^n:= R_{k+2(n-1)}\circ\dots \circ R_{k+2} \circ R_k.
\]
We also let $R_k^0$ be the identity.
The differential operator
\[
D:=\frac{1}{2\pi i} \frac{d}{dz}=q \frac{d}{dq}.
\]
satisfies the following relation
\begin{displaymath}
R_k=-4\pi D +k/y.
\end{displaymath}
The next lemma is often referred to as Bol's identity.

\begin{lemma}\label{bol}
Assuming the notation and hypotheses above, we have
\[
 D^{k-1}=\frac{1}{(-4\pi)^{k-1}}R_{2-k}^{k-1} .
\]
\end{lemma}

\begin{proof}
This is a special case of the identity (4.15) in \cite{LZ}.
\end{proof}

\section{Proof of Theorem~\ref{thm1}}\label{thm1proof}
By Lemma~\ref{bol},
we see that $D^{k-1}$ defines a linear map from
\begin{displaymath}
D^{k-1}: M^!_{2-k}(\Gamma_0(N),\chi)\longrightarrow
M^!_{k}(\Gamma_0(N),\chi).
\end{displaymath}
Theorem~\ref{thm1} asserts that this map may be extended to
harmonic weak Maass forms. Moreover, the theorem provides
a simple description of the images.

\begin{proof}[Proof of Theorem~\ref{thm1}]
Suppose that $k\geq 2$, and that $f\in H_{2-k}(\Gamma_0(N),\chi)$.
In view of Lemma \ref{bol}, it is clear that $D^{k-1} f$ has the
transformation behavior of a modular form of weight $k$.

We now show that $L_k D^{k-1} f =0$. This implies that $D^{k-1} f$ is holomorphic on $\H$.
By Lemma \ref{bol}, it suffices to show that  $L_k R_{2-k}^{k-1} f=0$.
Since $\Delta_{2-k} f=0$, it follows from \eqref{rev} by induction that
\[
\Delta_{k-2} R_{2-k}^{k-2} f = (2-k) R_{2-k}^{k-2} f.
\]
Using \eqref{deltalr}, we obtain
\begin{align*}
L_k R_{2-k}^{k-1} f= (L_k R_{k-2})R_{2-k}^{k-2} f
= (-\Delta_{k-2} -(k-2)) R_{2-k}^{k-2} f=0.
\end{align*}

Finally, the growth behavior of $f$ at the cusps implies that
$D^{k-1} f$ is meromorphic at the cusps. Therefore, $D^{k-1}$ indeed
extends to $H_{2-k}(\Gamma_0(N),\chi)$.

To complete the proof, we compute the Fourier expansion of
$D^{k-1}f$.
Assuming the notation in
\eqref{fourier}, a straightforward calculation gives
\begin{align*}
R_{2-k}^{k-2} f (z)= \sum_{n\gg -\infty} c_f^+(n) \Gamma(k-1,4\pi n y)(-y)^{2-k} e^{2\pi i n\bar z} + (k-2)!^2\sum_{n<0}
c_f^-(n)(-y)^{2-k}e^{2\pi i n\bar z}.
\end{align*}
Moreover,  $R_{2-k}^{k-1} f$ has the Fourier expansion
\begin{align*}
R_{2-k}^{k-1} f (z)=
\sum_{n\gg -\infty} c_f^+(n) (-4\pi n)^{k-1} q^n.
\end{align*}
In particular, we have
\begin{align*}
D^{k-1} f= D^{k-1} f^{+} =\sum_{n\gg -\infty} c_f^+(n) n^{k-1} q^n.
\end{align*}
The first two formulas follow from the Fourier expansion of $f$
and the differential equations
$\Delta_{k-2} R_{2-k}^{k-2} f = (2-k) R_{2-k}^{k-2} f$ and $\Delta_{k} R_{2-k}^{k-1} f = 0$.
The third formula is a consequence of the second and Lemma \ref{bol}.
\end{proof}

\begin{remark}
Note that $g:=y^{k-2}\overline{R_{2-k}^{k-2} f} $ is a harmonic weak Maass form of weight $2-k$ in the (slightly more general)
sense of Section 3 of \cite{BF}.
Moreover, $\xi_{2-k} g= y^{-k}\overline{L_{2-k} g} = R_{2-k}^{k-1} f$.
This can also be used to compute the Fourier expansions
in the proof of Theorem~\ref{thm1}.
\end{remark}

\section{The regularized inner product and the proof of Theorem~\ref{thm2}
}\label{innerprod}

Here we recall the regularized inner product, and we prove
Theorem~\ref{thm2}. We consider  slightly more
general situations, with earlier definitions modified in
the obvious way.

Let $k$ be an integer, and let $\Gamma$ be a subgroup of finite index of $\Gamma(1)=\Sl_2(\Z)$.
We define a regularized inner  product of $g\in M_k(\Gamma)$ and $h\in M^!_k(\Gamma)$ as follows.
For $T>0$ we denote by $\calF_T(\Gamma(1))$ the truncated fundamental domain
\[
\calF_T(\Gamma(1))=\{ z\in \H\ : \ \text{$|x|\leq 1/2$, $|z|\geq 1$, and $y\leq T$}\}
\]
for $\Gamma(1)$. Moreover, we define
the truncated fundamental domain for $\Gamma$ by
\[
\calF_T(\Gamma)=\bigcup_{\gamma\in \Gamma\bs \Gamma(1)}
\gamma\calF_T(\Gamma(1)).
\]
Following \cite{Bo1}, we define the regularized
inner product $( g,h)^{reg}$ as the constant term in the Laurent
expansion at $s=0$ of the meromorphic continuation in $s$ of the function
\[
\frac{1}{[\Gamma(1):\Gamma]} \lim_{T\to\infty}  \int_{\calF_T(\Gamma)} g(z)\overline{h(z)} y^{k-s}\,\frac{dx\,dy}{y^2}.
\]
Using the same argument as in Section 6 of \cite{Bo1}, it can be shown that $( g,h)^{reg}$ exists for any $g\in M_k(\Gamma)$ and $h\in M^!_k(\Gamma)$.
(It also exists for $g\in M_k(\Gamma)$ and $h\in H_k(\Gamma)$. But note that it does {\em not} exist in general if $g$ and $h$ are both weakly holomorphic with honest poles at the cusps.) It is clear, for cusp forms $g$ and $h$, that
the regularized inner product reduces to the classical Petersson inner product
$(g,h)$.

\begin{remark}
If  $h\in M^!_k(\Gamma)$ has vanishing constant term at every
cusp of $\Gamma$, then
\[
( g,h)^{reg} = \frac{1}{[\Gamma(1):\Gamma]} \lim_{T\to\infty}  \int_{\calF_T(\Gamma)} g(z)\overline{h(z)} y^k\,\frac{dx\,dy}{y^2}.
\]
\end{remark}

For the rest of this section we assume that $k\geq 2$.

\begin{theorem}
\label{inp} If $g\in M_k(\Gamma)$ and $f\in H_{2-k}(\Gamma)$, then
\[
(g,R_{2-k}^{k-1} f)^{reg} = \frac{(-1)^{k}}{[\Gamma(1):\Gamma]}\sum_{\kappa\in \Gamma\bs P^1(\Q)}
w_\kappa \cdot c_g(0,\kappa)\overline{c_f^+(0,\kappa)},
\]
where $c_g(0,\kappa)$ (resp. $c_f^+(0,\kappa)$) denotes the constant term of the Fourier expansion of $g$ (resp. $f$) at the cusp $\kappa\in P^1(\Q)$, and $w_\kappa$ is the width of the cusp $\kappa$.
\end{theorem}

\begin{proof}
For simplicity, we carry out the proof only in the special case $\Gamma=\Gamma(1)$.
The general case is completely analogous.
We put $H:=y^{k-2} \overline{R^{k-2}_{2-k} f}$. Then $h:=R^{k-1}_{2-k} f=y^{-k}\overline{L_{2-k} H}$.
Since the constant terms at all cusps of $h$ vanish, we have
\begin{align*}
(g,R_{2-k}^{k-1} f)^{reg} &= \lim_{T\to\infty}  \int_{\calF_T(\Gamma)} g(z)\overline{h(z)} y^k\,\frac{dx\,dy}{y^2}\\
&= \lim_{T\to\infty}  \int_{\calF_T(\Gamma)} g(z)(L_{2-k} H)\,\frac{dx\,dy}{y^2}\\
&=\lim_{T\to\infty}  \int_{\calF_T(\Gamma)} g(z)(\frac{\partial}{\partial \bar z} H)\,dz\,d\bar z\\
&=-\lim_{T\to\infty}  \int_{\calF_T(\Gamma)} (\bar \partial H)\wedge g(z)\,d z.
\end{align*}
Using the holomorphy of $g$, we obtain, by Stokes' theorem, the expression
\begin{align*}
(g,R_{2-k}^{k-1} f)^{reg} &=-\lim_{T\to\infty}  \int_{\calF_T(\Gamma)} d( H(z) g(z)\,d z)\\
&=-\lim_{T\to\infty}  \int_{\partial \calF_T(\Gamma)}  H(z) g(z)\,d z\\
&=\lim_{T\to\infty}  \int_{x=-1/2}^{1/2} H(x+iT) g(x+iT)\,d x.
\end{align*}
The integral over $x$ gives the constant term in the Fourier expansion of $H(x+iT) g(x+iT)$. It can be computed using the Fourier expansion
\begin{align*}
H(z)= (-1)^k\sum_{n\gg -\infty} \overline{c_f^+(n)} \Gamma(k-1,4\pi n y) e^{-2\pi i n z} + (-1)^k(k-2)!^2\sum_{n<0}
\overline{c_f^-(n)} e^{-2\pi i n z}
\end{align*}
of $H$ (see the proof of Theorem~\ref{thm1}) and the Fourier expansion of $g$.
It turns out that under the limit $T\to \infty$ only the contribution $(-1)^{k}
c_g(0)\overline{c_f^+(0)}$
coming from the product of the individual constant terms survives.
This concludes the proof.
\end{proof}

\begin{corollary}
If $g\in S_k(\Gamma)$, then $(g,R_{2-k}^{k-1} f)^{reg}=0$.
\end{corollary}

\begin{proof}
This is a direct consequence of Theorem~\ref{inp}.
\end{proof}

The next corollary implies Theorem~\ref{thm2}.

\begin{corollary}
The image of the map $D^{k-1}: H_{2-k}(\Gamma)\to M^!_{k}(\Gamma)$
is given by those $h\in M^!_k(\Gamma)$ which are orthogonal to cusp
forms and whose constant term at any cusp of $\Gamma$ vanishes.
\end{corollary}

\begin{proof}
If $f\in H_{2-k}(\Gamma)$, it follows from Theorem~\ref{inp} and
Theorem~\ref{thm1} that $D^{k-1}f$ satisfies the stated conditions.

Conversely, assume that $h\in M^!_k(\Gamma)$ is orthogonal to cusp forms and has vanishing constant term at any cusp
 of $\Gamma$.
According to Lemma 3.11 of \cite{BF}, we may chose $f\in
H_{2-k}(\Gamma)$ such that the principal parts of $D^{k-1}f$ and $h$
at the cusps agree up to the constant terms. Since the constant
terms of $h$ and $D^{k-1}f$ vanish, they trivially agree as well.
Consequently,
\[
h-D^{k-1}f\in S_{k}(\Gamma).
\]
In view of Theorem~\ref{inp} and the hypothesis on $h$, we find that $h- D^{k-1}f$ is orthogonal to cusp forms.
Hence it vanishes identically.
\end{proof}
\begin{remark}
It suffices to specialize $\Gamma = \Gamma_1(N)$ in the previous
Corollary in order to derive Theorem \ref{thm2}.
\end{remark}

\section{Proof of Theorems \ref{thm3}, \ref{thm:heckeVanishing},
and \ref{serrepadic}}\label{proof1.3}
Here we prove Theorem~\ref{thm3} by combining facts about
$\xi_{2-k}$, with Hecke theory and the theory of complex
multiplication. We first begin with an important proposition.

\begin{proposition}
\label{goodmaass} Let $g=\sum_{n=1}^{\infty}b(n)q^n\in
S_k(\Gamma_0(N),\overline{\chi})$ be a normalized newform with
integer weight $k\geq 2$, and let $F_g$ be the number field obtained
by adjoining the coefficients of $g$ to $\Q$. Then there is a
harmonic weak Maass form $f\in H_{2-k}(\Gamma_0(N), \chi)$ which
satisfies:
\begin{enumerate}
\item[(i)]
The principal part of $f$ at the cusp $\infty$ belongs to  $F_g[q^{-1}]$.
\item[(ii)]
The principal parts of $f$ at the other cusps of $\Gamma_0(N)$ are constant.
\item[(iii)]
We have that $\xi_{2-k}(f)=\|g\|^{-2} g$.
\end{enumerate}
\end{proposition}

\begin{proof}
Let $H_{2-k,\infty}(\Gamma_0(N),\chi)$ be the subspace of those
$f\in H_{2-k}(\Gamma_0(N),\chi)$ whose principal parts at the cusps
other than $\infty$ are constant. Note that
\begin{displaymath}
H_{2-k}(\Gamma_0(N),\chi)=H_{2-k,\infty}(\Gamma_0(N),\chi)
+M_{2-k}^!(\Gamma_0(N),\chi).
\end{displaymath}
Arguing as in Section 3 of \cite{BF}, the restriction of $\xi_{2-k}$
to $H_{2-k,\infty}(\Gamma_0(N),\chi)$
defines a surjective map to $S_k(\Gamma_0(N),\overline{\chi})$.
One now argues as in the proof of
Lemma 7.3 of \cite{BrO} using the pairing
$\{ g,f\}=(g,\xi_{2-k}(f))$, where
$f\in H_{2-k,\infty}(\Gamma_0(N),\chi)$
and $g\in S_k(\Gamma_0(N),\overline{\chi})$.
\end{proof}

\begin{remark}
The harmonic weak Maass form $f$ satisfying (i)--(iii) above is
unique up to the addition of a weakly holomorphic form  in
$M_{2-k}^!(\Gamma_0(N),\chi)$ with coefficients in $F_g$ and a pole
possibly at infinity and constant principal part at all other cusps.
\end{remark}

For completeness, here we briefly recall the notion of a newform
with complex multiplication (for example, see Chapter 12 of
\cite{iwaniec} or Section 1.2 of \cite{OnoCBMS}).
Let $D<0$ be the fundamental discriminant of an imaginary
quadratic field $K=\Q(\sqrt{D})$. Let $O_K$ be
the ring of integers of $K$, and let $\chi_K:=\leg{D}{\bullet}$
be the usual Kronecker character associated to $K$.
Let $k\geq 2$, and let
$c$ be a Hecke character of $K$ with exponent $k-1$ and
conductor $\mathfrak{f}_c$, a non-zero ideal of $O_K$.
By definition, this means that
\begin{displaymath}
c: I(\mathfrak{f}_c)\longrightarrow \C^{\times}
\end{displaymath}
is a homomorphism, where $I(\mathfrak{f}_c)$ denotes the group
of fractional ideals of $K$ prime to $\mathfrak{f}_c$.
In particular, this means that
\begin{displaymath}
c(\alpha O_K)=\alpha^{k-1}
\end{displaymath}
for $\alpha\in K^{\times}$ for which $\alpha\equiv 1\ {\text {\rm mod}}^{\times}
\mathfrak{f}_c$.
To $c$ we naturally associate a Dirichlet character $\omega_c$
defined, for every integer $n$ coprime to $\mathfrak{f}_c$, by
\begin{displaymath}
  \omega_c(n):=\frac{c(nO_K)}{n^{k-1}}.
\end{displaymath}
Given this data, we let
\begin{equation}\label{CMform}
\Phi_{K,c}(z):=\sum_{\mathfrak{a}} c(\mathfrak{a})q^{N(a)},
\end{equation}
where $\mathfrak{a}$ varies over the ideals of $O_K$ prime to
$\mathfrak{f}_c$, and where $N(\mathfrak{a})$ is the usual ideal norm.
It is well known that $\Phi_{K,c}(z)\in S_k(\Gamma_0(|D|\cdot
N(\mathfrak{f}_c)),\chi_K\cdot \omega_c)$ is a normalized newform.
These are newforms with {\it complex multiplication}.
By construction, if we let
\begin{displaymath}
\Phi_{K,c}(z)=\sum_{n=1}^{\infty}b(n)q^n,
\end{displaymath}
then
\begin{equation}\label{inert}
b(n)=0 \ \ \ {\text {\rm whenever}}\ \chi_K(n)=-1.
\end{equation}
This follows since every prime $p$ for which $\chi_K(p)=-1$ is inert.

\begin{proof}[Proof of Theorem \ref{thm3}.]
Suppose that $f$ is good for a CM form $g=\sum_{n=1}^{\infty}b(n)q^n$, and let
$D=D_g$ be the fundamental discriminant of the associated imaginary
quadratic field $K=\Q(\sqrt{D})$.
By (\ref{xiprop})
(correcting a typographical error in Lemma 3.1 of \cite{BF}), we then have that
\begin{displaymath}
\xi_{2-k}(f)=\|g\|^{-2}g=-(4\pi)^{k-1}\sum_{n=1}^{\infty}
\overline{c_f^{-}(-n)}n^{k-1}q^n.
\end{displaymath}
Since $g$ has complex multiplication, (\ref{inert}) implies
that
$c^-_f(n)=0$ when $\chi_{K}(-n)=-1$.
Because $D<0$, this means that
\begin{align}
\label{h1}
\text{$c^-_f(n)=0$ when $\chi_{K}(n)=1$}.
\end{align}

Let $M=ND$.
We write $\chi_0$ for the trivial character modulo $|D|$.
Since $D\mid N$, a standard argument shows that the sum of character twists
\[
u:=f\otimes\chi_0 + f\otimes \chi_{K}
\]
is in $H_{2-k}(\Gamma_0(M),\chi)$. The Fourier expansion of
$u=u^++u^-$ is given by
\begin{align*}
&u^+(z)=2\sum_{\substack{ n\gg -\infty\\ \chi_{K}(n)=1}} c_f^+(n) q^n,\\
&u^-(z)=2\sum_{\substack{n<0\\ \chi_{K}(n)=1}} c_f^-(n)\Gamma(k-1,4\pi|n|y)q^n.
\end{align*}
Consequently, by \eqref{h1}, the non-holomorphic part $u^-$
vanishes, and $u$ is actually weakly holomorphic.

We now claim that for any integer $b$, $f(z+b/D)$ has
principal parts at all cusps in $F_g(\zeta_{M})[q^{-1}]$.  To see
this, we let $\gamma\in \Gamma(1)$ and consider the cusp
$\gamma\infty$.  There exists a $\tilde\gamma\in \Gamma(1)$ and
$\alpha,\beta,\delta\in \Z$ such that
\[
\zxz{D}{b}{0}{D}\gamma =\tilde \gamma \zxz{\alpha}{\beta}{0}{\delta}.
\]
Hence, the Fourier expansion of $f(z+b/D)$ at the cusp $\gamma\infty$
is given by
\[
f\mid \tilde\gamma \mid \zxz{\alpha}{\beta}{0}{\delta}.
\]
By the assumption of $f$, it is holomorphic at the cusp $\infty$,
unless $\tilde \gamma\in \Gamma_0(N)$, in which case it is equal to
\[
f\mid \zxz{\alpha}{\beta}{0}{\delta}.
\]
Since $\delta\mid D^2\mid M$, the principal part at $\infty$ of this
modular form is contained in $F_g(\zeta_M)[q^{-1}]$, proving the claim.
This implies that  the twists $f\otimes \chi_0$,
$f\otimes \chi_D$, have
principal parts at all cusps in $F_g(\zeta_{M})[q^{-1}]$.
Therefore, the same is true for $u$.

Now we
recall the fact that
the action of $\Aut(\C/\Q(\zeta_N))$ commutes with the
action of $SL_2(\Z)$ on modular functions for $\Gamma(N)$
(for example, see
Theorem 6.6 in Chapter 6.2 and the diagram before
Remark 6.7 in Shimura's book \cite{shimurabook}).
Using the action of $\Aut(\C/F_g(\zeta_{M}))$ on
weakly holomorphic modular forms, we see that $u^\sigma$ has the
same properties for any $\sigma\in \Aut(\C/F_g(\zeta_M))$.
Moreover, $u^\sigma$ has the same principal parts as $u$ at all cusps.
Hence the
difference $u-u^\sigma$ is a holomorphic modular form which vanishes
at the cusp $\infty$. Since $2-k\leq 0$, this implies that
$u=u^\sigma$. Consequently, $u$ is defined over $F_g(\zeta_M)$. So
for all $n\in \Z$ with $\chi_{K}(n)=1$, we have that $c_f^+(n)\in
F_g(\zeta_M)$. In particular, $c_f^+(1)\in F_g(\zeta_M)$.

We now use the Hecke action on $f$ and $g$.
Let $T(m)$ be the $m$-th Hecke operator for $\Gamma_0(N)$.
Using the same argument as in Lemma 7.4 of \cite{BrO},  we have that
\[
f\mid_{2-k} T(m) = m^{1-k} b(m) f + f',
\]
where $f'\in M^!_{2-k}(\Gamma_0(N),\chi)$ is a weakly holomorphic
form with coefficients in $F_g$. In view of the formula for the
action of the Hecke operators on the Fourier expansion, we obtain
for any prime $p$ that
\[
c^+_f(pn) +{\chi(p)} p^{1-k} c^+_f(n/p) = p^{1-k} b(p) c^+_f(n)+ c^+_{f'}(n),
\]
where $c^+_{f'}(n)\in F_g$. Hence an inductive argument shows that
all coefficients $c^+_f(n)$ are contained in the extension
$F_g(c_f^+(1))$. This concludes the proof of the theorem since we
have already established that $c_f^{+}(1)$ is in $F_g(\zeta_M)$.
\end{proof}

The proof of Theorem \ref{thm:heckeVanishing} is similar to the
proof of Theorem \ref{thm3}, and so we give only a sketch
of the set-up.

\begin{proof}[Sketch proof of Theorem \ref{thm:heckeVanishing}]
If $p\nmid N$ is a prime, then for
every positive integer $m$ we have that
\begin{displaymath}
c_g(p)c_g(p^m)=c_g(p^{m+1})+\overline{\chi(p)}p^{k-1}c_g(p^{m-1}).
\end{displaymath}
Therefore, if $p\nmid N$ is a prime for which $c_g(p)=0$, then we have
that
\begin{displaymath}
c_g(p^{m+1})=-\overline{\chi(p)}p^{k-1}c_g(p^{m-1}),
\end{displaymath}
which in turn implies that
\begin{displaymath}
c_g(p^{m})=\begin{cases} \left(-\overline{\chi(p)}p^{k-1}\right)^{\frac{m}{2}}
\ \ \ \ \ &{\text {\rm if}}\ m\ {\text {\rm is even}},\\
0 \ \ \ \ \ &{\text {\rm otherwise.}}
\end{cases}
\end{displaymath}
Therefore, by arguing with the usual $U(p), V(p), U(p^2)$ and $V(p^2)$
operators, we can obtain a harmonic weak Maass form whose Fourier
coefficients are supported on terms whose exponents $n$ have the property
that $p$ exactly divides $n$. By the multiplicativity of the Fourier
coefficients of newforms, it then follows by the observation
above that the non-holomorphic part of this form is identically zero.
In other words, this particular harmonic weak Maass form is a weakly
holomorphic modular form with suitable principal parts at cusps.
The proof now
follows {\it mutatis mutandis} as in the proof of Theorem~\ref{thm3}.
\end{proof}

\subsection{Proof of Theorem~\ref{serrepadic}}\label{padic}
Theorem~\ref{serrepadic} follows easily from the seminal work of
Serre \cite{Serre1, Serre2} on $p$-adic modular forms.
We now apply his works to prove the theorem.

By Theorem~\ref{thm1}, we have that
\begin{displaymath}
\sum_{n\gg -\infty} c_f^{+}(n)n^{k-1}q^n
\end{displaymath}
is a weight $k$ weakly holomorphic modular form on $\Gamma_0(p^t)$
with rational coefficients. By standard facts involving the
$U(p)$-operator
\begin{displaymath}
\left(\sum a(n)q^n\right)\ | \ U(p):=\sum a(np)q^n,
\end{displaymath}
the integer $a$ has the property that
\begin{equation}
f^{*}(z):=\sum_{n=0}^{\infty} c_f^{+}(p^an) n^{k-1}q^n
\end{equation}
is a weight $k$ weakly holomorphic modular form on
$\Gamma_0(p^{t^{*}})$, where $t^*=1$ if $t=0$, and is $t$ otherwise.
This modular form has trivial principal part at the cusp infinity.
Therefore, we may apply a theorem of Serre (see Th.~5.4 of
\cite{Serre2}), and the conclusion is that $f^{*}(z)$ is a $p$-adic
modular form on $\SL_2(\Z)$ of weight $k$. This proves (1).

Claim (2) follows from the definition of a $p$-adic modular form.
Indeed, $p$-adic modular forms are $p$-adic limits of the Fourier
expansions of classical holomorphic modular forms, and these
forms, by a theorem of Serre,
have the property that almost all of their coefficients
are multiples of any fixed power of $p$ (see Th. 4.7 of
\cite{Serre2}).
This implies (2).

Claim (3) is a consequence of the fact that the $U(p)$-operator acts
locally nilpotently on certain $p$-adic modular forms. In this
situation, Serre proved that the constant term of a $p$-adic modular
form is a limit, in the $p$-adic sense, of certain Fourier
coefficients. This result (see Th. 7 and the following
remark in \cite{Serre1}) implies (3).

\section{Poincar\'e series and Theorems~\ref{thm1}, \ref{thm2}, and
\ref{thm3}
}\label{Poincare}

Here we consider natural examples of the results of both
Theorems~\ref{thm1}, \ref{thm2} and \ref{thm3}. Our results depend
on the explicit Fourier expansions
of two classes of Poincar\'e series (for example, see \cite{
BOPNAS, Br, iwaniec}).

\subsection{Definitions and Fourier expansions}
For $A=\kabcd\in \SL_2(\Z)$, define $j(A,z)$ by
\begin{equation}
j(A,z):=(cz+d).
\end{equation}
As usual, for such $A$ and functions $f:\H\to
\C$, we let
\begin{equation}\label{slash}
(f\mid_k A )(z):= j(A,z)^{-k} f(A z).
\end{equation}
Let $m$ be an integer, and let $\varphi_m:\R^{+}\to \C$ be a
function which satisfies $\varphi_m(y)=O(y^\alpha)$, as $y\to 0$,
for some $\alpha\in \R$. If $e(\alpha):=e^{2\pi i \alpha}$ as
before, then let
\begin{equation}
\varphi^{*}_m(z):=\varphi_m(y)e(mx).
\end{equation}
Such functions are fixed by the translations $\Gamma_\infty:=\{ \pm
\kzxz{1}{n}{0}{1}\ :\ n\in \Z\}$.

Given  this data, for
integers $N\geq 1$, we define
the generic Poincar\'e series
\begin{align}
\PP(m,k,\varphi_m,N;z):= \sum_{A\in\Gamma_\infty \bs \Gamma_0(N)}
(\varphi^{*}_m \mid_k A)(z).
\end{align}
We shall be interested in two families of such series.

The first family is classical.
We let
\begin{equation}\label{classicalpoincare}
P(m,k,N;z)=q^m+\sum_{n=1}^{\infty}a(m,k,N;n)q^n:=\PP(m,k,e(imy),N;z).
\end{equation}
These series are modular, and their Fourier expansions
are given in terms of the $I$-Bessel and $J$-Bessel functions, and
the Kloosterman sums
\begin{equation}
K(m,n,c):=
\sum_{v(c)^{\times}} e\left(\frac{m\overline
v+nv}{c}\right).
\end{equation}
Here $v$ runs through the primitive residue classes modulo $c$,
and $v\overline v\equiv 1\pmod{c}.$
The following is well known (for example, see \cite{iwaniec,
Petersson}).

\begin{proposition}\label{poinc1prop}
If $k\in 2\N$, and $m, N\geq 1$, then the following are true.

\smallskip
\noindent 1) We have that $P(m,k,N;z)\in S_k(\Gamma_0(N))$, and for
positive integers $n$ we have
\begin{displaymath}
a(m,k,N;n)=2\pi (-1)^{\frac{k}{2}}\left(\frac{n}{m}\right)^{\frac{k-1}{2}}
\cdot
\sum_{\substack{c>0\\c\equiv 0\pmod{N}}} \frac{K(m,n,c)}{c}\cdot J_{k-1}
\left(\frac{4\pi \sqrt{mn}}{c}\right).
\end{displaymath}

\smallskip
\noindent 2) We have that $P(-m,k,N;z)\in M_k^{!}(\Gamma_0(N))$, and
for positive integers $n$ we have
\begin{displaymath}
a(-m,k,N;n)=2\pi(-1)^{\frac{k}{2}}\left(\frac{n}{m}\right)^{\frac{k-1}{2}}
\cdot \sum_{\substack{c>0\\ c\equiv 0\pmod{N}}}\frac{K(-m,n,c)}{c}\cdot
I_{k-1}\left(\frac{4\pi \sqrt{|mn|}}{c}\right).
\end{displaymath}
\end{proposition}

Now we recall the second family of Poincar\'e series, the
Maass-Poincar\'e series of Hejhal (see \cite{Hejhal}).
Let
$M_{\nu,\,\mu}(z)$
be the usual $M$-Whittaker
function.
For complex $s$,
let
$$\calM_s(y):=
|y|^{-\frac{k}{2}} M_{\frac{k}{2}\sgn(y),\,s-\frac{1}{2}}(|y|),
$$
and for $m\geq 1$ let
$\varphi_{-m}(z):=\calM_{1-\frac{k}{2}}(-4\pi m y)$.
For $k\in 2\N$ and integers $N\geq 1$, we let
\begin{equation}
Q(-m,k,N;z):=\frac{1}{(k-1)!}\cdot \PP(-m,2-k,\varphi_{-m},N;z).
\end{equation}
To determine the Fourier expansions of these series, we also require
the incomplete Gamma-function $\Gamma(a,x)$. We have the following
proposition (for example, see \cite{BOPNAS, Br, Hejhal}).

\begin{proposition}\label{poinc2prop}
If $k\in 2\N$, and $m, N\geq 1$, then $Q(-m,k,N;z)\in
H_{2-k}(\Gamma_0(N))$, and has a Fourier expansion of the form
\begin{displaymath}
Q(-m,k,N;z)= Q^{+}(-m,k,N;z) + Q^{-}(-m,k,N;z),
\end{displaymath}
where
\begin{displaymath}
Q^{-}(-m,k,N;z) =-\frac{\Gamma(k-1,4\pi m y)}{(k-2)!}q^{-m}+
\sum_{n<0}b(-m,k,N;n) \cdot\Gamma(k-1,4\pi |n|y) q^n,
\end{displaymath}
and where for negative integers $n$ we have
\begin{displaymath}
b(-m,k,N;n) =  - \frac{2\pi (-1)^{\frac{k}{2}}}{(k-2)!} \cdot \left
| \frac{m}{n}\right|^{\frac{k-1}{2}} \sum_{\substack{c>0\\ c\equiv
0\pmod{N}}} \frac{K(-m,n,c)}{c}\cdot J_{k-1}\left( \frac{4\pi \sqrt{
|mn|}}{c}\right),
\end{displaymath}
and
\begin{displaymath}
Q^{+}(-m, k, N;z)=q^{-m}+\sum_{n=0}^{\infty}b(-m,k,N;n)q^n,
\end{displaymath}
where
\begin{displaymath}
b(-m,k,N;0)=-\frac{2^k\pi^k (-1)^{\frac{k}{2}}m^{k-1}}{(k-1)!}\cdot
\sum_{\substack{c>0\\ c\equiv 0\pmod{N}}}\frac{K(-m,0,c)}{c^k},
\end{displaymath}
and where for positive integers $n$ we have
\begin{displaymath}
b(-m,k,N;n)=-2\pi (-1)^{\frac{k}{2}}
\cdot \sum_{\substack{c>0\\ c\equiv 0\pmod{N}}}
\left(\frac{m}{n}\right)^{\frac{k-1}{2}}
\frac{K(-m,n,c)}{c}\cdot I_{k-1}\left(\frac{4\pi \sqrt{|mn|}}{c}\right).
\end{displaymath}
\end{proposition}

\begin{remark} Obviously, we have that  $Q^{+}(-m,k,N;z)$
(resp. $Q^{-}(-m,k,N;z)$)
is the holomorphic part (resp. non-holomorphic part) of the weak Maass form
$Q(-m,k,N;z)$.

\begin{remark} Propositions~\ref{poinc1prop} and \ref{poinc2prop} are well known for
 $2<k\in 2\N$. That they hold for $k=2$ follows by arguing by
analytic continuation in $k$ with Fourier expansions.
\end{remark}

\end{remark}
\subsection{Poincar\'e series in the context of Theorems~\ref{thm1} and ~\ref{thm2}}

If $k\geq 2$ is even, and $m, N\geq 1$, then
Theorem~\ref{thm1}, and Propositions ~\ref{poinc1prop} and
~\ref{poinc2prop}, imply that
\begin{displaymath}
\begin{split}
&D^{k-1}Q(-m,k,N;z)=D^{k-1}Q^{+}(-m,k,N;z)=\\
&\ \ =(-m)^{k-1}q^{-m}-2\pi (-1)^{\frac{k}{2}}\cdot \sum_{n=1}^{\infty}
\sum_{\substack{c>0\\ c\equiv 0\pmod{N}}} (mn)^{\frac{k-1}{2}}\frac{K(-m,n,c)}{c}\cdot
I_{k-1}\left(\frac{4\pi \sqrt{|mn|}}{c}\right)\\
&\ \ =-m^{k-1}P(-m,k,N;z).
\end{split}
\end{displaymath}
In other words, we have
\begin{equation}\label{corr}
D^{k-1}Q(-m,k,N;z)=-m^{k-1}P(-m,k,N;z).
\end{equation}
By (\ref{xiprop}), we also have that
\begin{displaymath}
\xi_{2-k}(Q(-m,k,N;z))= \frac{( 4\pi)^{k-1}
m^{k-1}}{(k-2)!}\cdot P(m,k,N;z).
\end{displaymath}
For one dimensional spaces $S_k(\Gamma_0(N))$, say generated by
a newform $g$, this last relation
relates Maass-Poincar\'e series to newforms $g$.
This follows from standard facts about Petersson inner products
and Poincar\'e series (for example, see Chapter 3 of \cite{iwaniec}).

\section{A good example}
Correspondence (\ref{corr}) can be useful for computing Fourier
coefficients of holomorphic parts of certain harmonic weak Maass forms.
Here we consider example (\ref{goodpoincare}) from the introduction.
Using Proposition~\ref{poinc2prop}, we summed the first
150 terms to obtain
\begin{displaymath}
\begin{split}-D^3Q(-1,&4,9;z)
=P(-1,4,9;z)\\
&\ \ \ \ \ \ \ \ \ \ \sim q^{-1}+1.9999 q^2 -48.9999q^5+47.9999q^8+770.9999
q^{11}+\cdots.
\end{split}
\end{displaymath}
On the other hand, we have the weight 4 weakly holomorphic modular
form
\begin{displaymath}
m(z):=\left(\frac{\eta(z)^3}{\eta(9z)^3}+3\right)^2\cdot \eta(3z)^8=
q^{-1}+2q^2-49q^5+48q^8+771q^{11}-\cdots,
\end{displaymath}
where $\eta(z):=q^{\frac{1}{24}}\prod_{n=1}^{\infty}(1-q^n)$ is
Dedekind's eta-function. These two modular forms are equal, and so
the coefficients are obviously rational (in fact, integral). To
deduce this, one may use the ``circle method'' to get asymptotics
for the coefficients of $m(z)$ (for example, see the detailed
discussion in \cite{Rademacher}). 
The circle method gives the same asymptotic expressions for the coefficients
of $m(z)$ and $P(-1,4,9;z)$. 
This follows from the fact that they have the same principal parts
at cusps. The circle method
then shows that these approximations for the $n$th coefficients
of $m(z)$
and $P(-1,4,9;z)$  agree up to a power of $n$. This follows from a standard
calculation involving bounds for Kloosterman sums and the asymptotic
properties of $I$-Bessel functions. 
Therefore, it follows that
$P(-1,4,9;z)-m(z)$ 
is a holomorphic modular form. This
form vanishes at all cusps, and so it must be a cusp form.
The space
$S_4(\Gamma_0(9))$ is one dimensional, and is spanned by the CM form
\begin{displaymath}
g(z)=\eta(3z)^8=q-8q^4+20q^7-70q^{13}+64q^{16}+56q^{19}+\cdots.
\end{displaymath}
The non-zero coefficients of this cusp form
have exponents $n$ in the arithmetic
progression $n\equiv 1\pmod{3}$. However, since $K(-1,n,9c)=0$ for
integers $n\equiv 1\pmod 3$, Proposition~\ref{poinc1prop} (2) then
implies that $P(-1,4,9;z)-m(z)$ is identically zero.
Combining these facts, we find that the coefficients of
$Q^{+}(-1,4,9;z)$ are rational, and its first few terms are
\begin{displaymath}
Q^{+}(-1,4,9;z)=q^{-1}-\frac{1}{4}q^2+\frac{49}{125}q^5-\frac{48}{512}q^8-\frac
{771}{1331}q^{11}+\cdots.
\end{displaymath}

In the context of Theorem~\ref{thm3}, the algebraicity of
$Q^{+}(-1,4,9;z)$ follows from the fact that $Q(-1,4,9;z)$ is good
for the CM form $g(z)$. To see this one observes that
$\xi_{-2}(Q(-1,4,9;z))=\|g\|^{-2}g$. This follows from Petersson's
theory, combined with the fact that $S_4(\Gamma_0(9))$ is one
dimensional, and is spanned by both $g(z)$ and $P(1,4,9;z)$.

Theorem~\ref{serrepadic} is also easily described in this example.
We have that $a=1$ in Theorem~\ref{serrepadic}, and so
the $q$-series $f^{*}$
in Theorem~\ref{serrepadic} (1) is $-m(z)\ | U(3)$.
Theorem~\ref{serrepadic} (2)  then implies that almost every
coefficient of $m(z)\ | U(3)$ is a multiple of any fixed power of 3.
Theorem~\ref{serrepadic} (3) follows trivially in this case
since the coefficients of $m(z)$ for exponents which are powers
of $3$ obviously vanish.

\end{document}